\documentclass{article}
\usepackage{amsmath}
\usepackage{amssymb}
\newtheorem{thm}{Theorem}[section]

\newtheorem{cor}[thm]{Corollary}
\title{Poincar\'e series of some hypergraph algebras}
\author{E.~Emtander, Stockholm university\\
erice@math.su.se
\and
R.~Fr\"oberg, Stockholm university\\
ralff@math.su.se
\and
F.~Mohammadi, Amirkabir University, Tehran\\
f\_mohammadi@aut.ac.ir
\and
S.~Moradi Amirkabir University, Tehran\\
s\_moradi@aut.ac.ir}

\begin{document}
\maketitle
\begin{abstract}
A hypergraph $H=(V,E)$, where $V=\{ x_1,\ldots,x_n\}$ and $E\subseteq 2^V$
defines a hypergraph algebra $R_H=k[x_1,\ldots,x_n]/(x_{i_1}\cdots x_{i_k};
\{ i_1,\ldots,i_k\}\in E)$. All our hypergraphs are $d$-uniform, i.e., 
$|e_i|=d$ for all $e_i\in E$. We determine the Poincar\'e series
$P_{R_H}(t)=\sum_{i=1}^\infty\dim_k{\rm Tor}_i^{R_H}(k,k)t^i$ for some
hypergraphs generalizing lines, cycles, and stars. We finish by calculating
the graded Betti numbers and the Poincar\'e series of the graph algebra
of the wheel graph.
\end{abstract}
\section{Introduction}
A line is a graph $L_n=(V,E)$, where 
$$V=\{ x_1,\ldots,x_{n+1}\}\mbox{ and }E=
\{(x_1,x_2),\ldots,(x_n,x_{n+1})\},$$ a cycle a graph $C_n=
(V,E)$, where 
$$V=\{ x_1,\ldots,x_n\}\mbox{ and }E=
\{(x_1,x_2),\ldots,(x_{n-1},x_n),(x_n,x_1)\},$$ and a star a graph
$S_n=(V,E)$, where 
$$V=\{ x_1,\ldots,x_{n+1}\}\mbox{ and }E=
\{(x_1,x_2),\ldots,(x_1,x_{n+1})\}.$$
In \cite[Chapter 7]{ja} the Betti numbers of their graph algebras, 
$$k[x_1,\ldots,x_{n+1}]/(x_1x_2,x_2x_3,\ldots,x_nx_{n+1}),$$
$$k[x_1,\ldots,x_n]/(x_1x_2,x_2x_3,\ldots,x_{n-1}x_n,x_nx_1),$$ 
and
$$k[x_1,\ldots,x_{n+1}]/(x_1x_2,x_1x_3,\ldots,x_1x_{n+1})$$
are determined. This is generalized to certain
``hyperlines'', ``hypercycles'', and ``hyperstars'' in \cite{em-mo-mo}. 
Here a hyperline is hypergraph with $nd-(n-1)\alpha$ vertices 
and $n$ edges $e_1,\ldots,e_n$, where all edges $e_1,\ldots,e_n$
have size $d$, and $e_i\cap e_j\ne\emptyset$ and has size $\alpha$ if and only 
if $|i-j|=1$, a hypercycle is hypergraph with $n(d-\alpha)$
vertices and $n$ edges $e_1,\ldots,e_n$, 
where all edges have size $d$, and $e_i\cap e_j\ne\emptyset$
and has size $\alpha$ if and only if $|i-j|=1\pmod n$, and the hyperstar
is hypergraph with $n(d-\alpha)$
vertices and $n$ edges $e_1,\ldots,e_n$, 
where all edges have size $d$, and for all $i,j$
$|e_i\cap e_j|=|\cap_{i=1}^ne_i|=\alpha>0$.
We denote the line hypergraph and 
its algebra with $L_n^{d,\alpha}$, the cycle hypergraph 
and its algebra with $C_n^{d,\alpha}$, and the star hypergraph and its
algebra $S_n^{d,\alpha}$. Their Betti numbers were determined
in \cite[Chapter 3]{em-mo-mo} (in the first two cases with the restriction
$2\alpha\le d$). In this paper we will 
determine the Poincar\'e series for the same algebras.
The Poincar\'e series of a graded $k$-algebra $R=k[x_1,\ldots,x_n]/I$ is 
$P_R(t)=\sum_{i=1}^\infty\dim_k{\rm Tor}_i^R(k,k)t^i$. \cite{gu-le} is an 
excellent source for results on Poincar\'e series.

\section{Hypercycles and hyperlines when $d=2\alpha$}
We start with the case $d=2\alpha$. If 
$e_i=\{ v_{i1},\ldots,v_{i\alpha},v'_{i1},\ldots,v'_{i\alpha}\}$, 
where $\{ v_{ij}'\}\in e_{i+1}$, we start by factoring out
all $v_{ik}-v_{il}$ and $v'_{ik}-v'_{il}$. This is a linear 
regular sequence of length $(n+1)(\alpha-1)$ for the hyperline
and of length $n(\alpha-1)$ for the hypercycle. The results are
$$L'_{n,a}=k[x_1,\ldots,x_{n+1}]/
(x_1^\alpha x_2^\alpha,x_2^\alpha x_3^\alpha,\ldots,x_n^\alpha x_{n+1}^\alpha)$$
and
$$C'_{n,a}=k[x_1,\ldots,x_n]/
(x_1^\alpha x_2^\alpha,x_2^\alpha x_3^\alpha,\ldots,
x_{n-1}^\alpha x_n^\alpha,x_n^\alpha x_1^\alpha).$$
Then 
$$P_{L_n^{2a,a}}(t)=(1+t)^{(n+1)(\alpha-1)}P_{L'_{n,a}}(t)$$
and 
$$P_{C_n^{2a,a}}(t)=(1+t)^{n(\alpha-1)}P_{C'_{n,a}}(t).$$ 
\cite[Theorem 3.4.2 (ii)]{gu-le}.
Now $L_n^{2\alpha,\alpha}$ and $C_n^{2\alpha,\alpha}$
have the same Poincar\'e series as the graph algebras
$$L_n=L_n^{2,1}=k[x_1,\ldots,x_{n+1}]/(x_1x_2,x_2x_3,\ldots,x_nx_{n+1})$$ 
and
$$C_n=C_n^{2,1}=k[x_1,\ldots,x_n]/(x_1x_2,x_2x_3,\ldots,x_{n-1}x_n,x_nx_1)$$
respectively.

For a graded $k$-algebra $\oplus_{i=0}^\infty R_i$, the Hilbert series of $R$
is defines as $H_R(t)=\sum_{i=0}^\infty\dim_k(R_i)t^i$. The exact sequences
$$0\longrightarrow (x_{n+1})\longrightarrow L_n\buildrel{x_{n+1}\cdot}
\over{\longrightarrow}
L_n\longrightarrow L_n/(x_{n+1})\longrightarrow 0$$ 
and
$$0\longrightarrow(x_{n+1})\longrightarrow L_n\longrightarrow L_n/(x_{n+1})
\longrightarrow 0$$
and 
$L_n/(x_{n+1})\simeq L_{n-1}$ and $(x_{n+1})\simeq L_{n-2}\otimes k[x]$
gives
$$H_{L_n}(t)=H_{L_{n-1}}(t)+\frac{t}{1-t}H_{L_{n-2}}(t).$$
The exact sequences
$$0\longrightarrow(x_1,x_{n-1})\longrightarrow C_n\buildrel{x_n\cdot}
\over{\longrightarrow}
C_n\longrightarrow L_{n-2}\longrightarrow 0$$
and
$$0\longrightarrow(x_1,x_{n-1})\longrightarrow C_n
\longrightarrow C_n/(x_1,x_{n-1})\longrightarrow 0$$
and
$C_n/(x_1,x_{n-1})\simeq L_{n-4}\otimes k[x]$ gives
$$H_{C_n}(t)=H_{L_{n-2}}(t)-{t\over(1-t)}H_{L_{n-4}}(t).$$
Now $C_n$ and $L_n$ are (as all graph algebras) Koszul algebras \cite[Corollary 2]{fr1}, so 
$P_{C_n}(t)=1/H_{C_n}(-t)$ and $P_{L_n}(t)=1/H_{L_n}(-t)$.
Since $L_0=k[x_1]$ and $L_1=k[x_1,x_2]/(x_1x_2)$, we have 
$H_{L_0}(t)=1/(1-t)$ and 
$H_{L_1}(t)=(1+t)/(1-t)$. We give the first Hilbert series:

\noindent
$H_{L_2}(t)=(1+t-t^2)/(1-t)^2$, $H_{L_3}(t)=(1+2t)/(1-t)^2$, 

\noindent
$H_{L_4}(t)=(1+2t-t^2-t^3)/(1-t)^3$,
$H_{L_5}(t)=(1+3t+t^2-t^3)/(1-t)^3$, 

\noindent
$H_{C_3}(t)=(1+2t)/(1-t)$, 
$H_{C_4}(t)=(1+2t-t^2)/(1-t)^2$,

\noindent
$H_{C_5}(t)=(1+3t+t^2)/(1-t)^3$, $H_{C_6}(t)=(1+3t-2t^3)/(1-t)^3$.

\medskip
\noindent{\bf Remark} We note that it is probably hard to get one
formula for $H_{L_n}(t)$ for all $n$. An indication is that we get 
the Fibonacci numbers from $H_{L_n}(t)$. For $t=1/2$
we get $H_{L_n}(1/2)=F_{n+2}$, the $(n+2)$th Fibonacci number if $F_0=F_1=1$. 

\medskip
Thus we get 

\noindent
$P_{L_2}(t)=(1+t)^2/(1-t-t^2)$, $P_{L_3}(t)=(1+t)^2/(1-2t)$,

\noindent 
$P_{L_4}(t)=(1+t)^3/(1-2t-t^2+t^3)$, $P_{L_5}(t)=(1+t)^3/(1-3t+t^2+t^3)$,

\noindent
$P_{C_3}(t)=(1+t)/(1-2t)$, 
$P_{C_4}(t)=(1+t)^2/(1-2t-t^2)$,

\noindent
$P_{C_5}(t)=(1+t)^2/(1-3t+t^2)$, $P_{C_6}(t)=(1+t)^3/(1-3t+2t^3)$.

\medskip
We collect the results in
\begin{thm}
The Poincar\'e series of $L_n$ and $C_n$ satisfy the recursion formulas
$$P_{L_n}(t)=\frac{(1+t)P_{L_{n-1}}(t)P_{L_{n-2}}(t)}{(1+t)P_{L_{n-2}}(t)-
tP_{L_{n-1}}(t)}$$
where $P_{L_0}(t)=1+t$ and $P_{L_1}(t)=(1+t)/(1-t)$ and
$$P_{C_n}(t)=\frac{(1+t)P_{L_{n-2}}(t)P_{L_{n-4}}(t)}{P_{L_{n-2}}(t)+
(1+t)P_{L_{n-4}}(t)}.$$
Furthermore
$$P_{L_n^{2\alpha,\alpha}}(t)=(1+t)^{(n+1)(\alpha-1)}P_{L_n}(t)$$ 
and 
$$P_{C_n^{2\alpha,\alpha}}(t)=(1+t)^{n(\alpha-1)}P_{C_n}(t).$$ 
\end{thm}

\section{Hypercycles and hyperlines when $2\alpha<d$}
Next we turn to the case $2\alpha<d$. Now each edge has a free vertex, i.e. a 
vertex which does not belong to any other edge.  
Then the Taylor resolution is minimal. 
In this case there is a formula for the Poincar\'e series in terms of
the graded homology of the Koszul complex \cite[Corollary to Proposition 2]{fr2}. Let $R$ be a monomial
ring for which the Taylor resolution is minimal.
Then the homology of the Koszul complex $H(K_R)$ is of the form
$H(K_R)=k[u_1,\ldots,u_N]/I$, where $I$ is generated by a set of monomials
of degree 2. Define a bigrading induced by $\deg(u_i)=(1,|u_i|)$, where 
$|u_i|$ is the homological degree. Then $P_R(t)=(1+t)^e/H_R(-t,t)$, where
$e$ is the embedding dimension and 
$H_R(x,y)$ is the bigraded Hilbert series of $H(K_R)$, see \cite{fr2}.

We begin with the hypercycle.
The homology of the Koszul complex is generated by $\{ z_I\}$, where 
$I=\{ i,i+1,\ldots,j\}$ corresponds to a
path $\{ e_i,e_{i+1},\ldots,e_j\}$ in $C_n^{d,\alpha}$ 
(indices counted $\pmod n$).
Thus there are $n$ generators in all homological degrees $<n$ and one
generator in homological degree $n$. We have $z_Iz_J=0$ if 
$I\cap J\ne\emptyset$.
Thus the surviving monomials are of the form $m=z_{I_1}\cdots z_{I_r}$, where
$I_i\cap I_j=\emptyset$ if $i\ne j$. The bidegree of $m$ is 
$(r,\sum_{j=1}^r|I_j|)$. Let $\sum_{j=1}^r|I_j|=i$. Then $m$ lies in 
$H(K)_{i,di-(i-r)\alpha}$. 
The graded Betti numbers are determined in \cite[Chapter 3]{em-mo-mo}.
The nonzero Betti numbers are $\beta_{i,di-(i-r)\alpha}=
\frac{n}{r}{i-1\choose r-1}{n-i-1\choose r-1}$ if $1\le r\le i<n$ and 
$\beta_{n,n(d-\alpha)}=1$. (As usual ${a\choose b}=0$ if $b>a$.) This gives the
Poincar\'e series.

Next we consider the hyperline. The homology of the Koszul complex is 
generated by $\{ z_I\}$, where $I=\{ i,i+1,\ldots,j\}$ corresponds to a
path $\{ e_i,e_{i+1},\ldots,e_j\}$ in $L(n,d,\alpha)$. Thus there are $n+1-i$
generators of homological degree $i$. We have $z_Iz_J=0$ if 
$I\cap J\ne\emptyset$. The graded Betti numbers are determined in \cite[Chapter 3]{em-mo-mo}.
The nonzero Betti numbers are $\beta_{i,di-(i-r)\alpha}=
{i-1\choose r-1}{n-i+1\choose r}$ if $1\le r\le i\le n$. The same reasoning
as above gives the Poincar\'e series. We state the results in a theorem.
\begin{thm}
If $2\alpha<d$, then
$$P_{C_n}(t)=\frac{(1+t)^{n(d-\alpha)}}{1+\sum_{1\le r\le i<n}(-1)^r\frac{n}{r}
{i-1\choose r-1}{n-i-1\choose r-1}t^{i+r}-t^{n+1}},$$
and
$$P_{L_n}(t)=\frac{(1+t)^{n(d-\alpha)+\alpha}}{1+\sum_{1\le r\le i\le n}
(-1)^r{i-1\choose r-1}{n-i+1\choose r}t^{i+r}}.$$
\end{thm}

\section{The hyperstar}
We conclude with a hypergraph generalizing the star graph. Suppose $|e_i|=d$
for all $i$, $1\le i\le n$, and that if $i\ne j$, then $|e_i\cap e_j|=
|\cap_{i=1}^ne_i|=\alpha<d$. Then the ideal is of the form $m(m_1,\ldots,m_n)$,
where $m$ is a monomial of degree $\alpha$. 
Then the hypergraph ring $S_n^{d,\alpha}$ is Golod
\cite[Theorem4.3.2]{gu-le}. This means that 

\begin{thm}
$$P_{S_n^{d,\alpha}}(t)=(1+t)^{|V|}/(1-\sum\beta_it^{i+1})=
(1+t)^{n(d-\alpha)+\alpha}/(1-\sum{n\choose i}t^{i+1}).$$
\end{thm}

\section{The wheel graph}
Finally we consider the wheel graph $W_n$, which is $C_n$ with an extra vertex
(the center) which is connected to all vertices in $C_n$. We let $W_n$ also
denote the graph algebra $k[x_0,\ldots,x_n]/(x_1x_2,x_2x_3,\ldots,x_nx_1,
x_0x_1,\ldots,x_0x_n)$.

\begin{thm}
Let $W_n$ be a wheel graph on $n+1$ vertices. Then the Betti numbers
of $W_n$ are as follows:\\

$(i)$ If $j>i+1$, then
$\beta_{i,j}(k[\Delta_{W_n}])=\beta_{i,j}(C_n)+\beta_{i-1,j-1}(C_n)$.

$(ii)$   If $j=i+1$, then
$\beta_{i,i+1}(W_n)=\beta_{i,i+1}(C_n)+\beta_{i-1,i}(C_n)
+{n\choose i}$.
\end{thm}
{\it Proof.}
{Assume that $V(W_n)=\{x_0,x_1,\ldots,x_n\}$ and $C_n=W_n\setminus
\{x_0\}$. It is easy to see that $\Delta_{W_n}=\Delta_{C_n}\cup
\{x_0\}$, where $\Delta_{W_n}$ and $\Delta_{C_n}$ are the
independence complexes of $W_n$ and $C_n$. It implies that for any
$i\geq 1$, $H_i(\Delta_{W_n})=H_i(\Delta_{C_n})$. Thus, if $j>i+1$, from
Hochster's formula (\cite[Theorem 5.5.1]{b-h}) 
and the observation above one has the result. Now
assume that $j=i+1$. Then
\\
$\beta_{i,i+1}(W_n)=\sum_{S\subseteq V(W_n),|S|=i+1}
\dim(\widetilde{H}_0(\Delta_S))=\sum _{S\subseteq V(C_n),|S|=i+1}
\dim(\widetilde{H}_0(\Delta_S))\\+\sum_{S\subseteq V(W_n),S=S'\cup
\{x_0\} } \dim(\widetilde{H}_0(\Delta_S))$. For any $S\subseteq
V(W_n)$ and $S_0\subseteq V(C_n)$, let $r_S$ and $r'_{S_0}$ denotes
the number of connected components of $\Delta_S$ in $V(W_n)$ and
$\Delta_{S_0}$ in $V(C_n)$ respectively. Then we have
$\sum_{S\subseteq V(W_n),S=S_0\cup \{x_0\} }
\dim(\widetilde{H}_0(\Delta_S))=\sum_{S\subseteq V(W_n),S=S_0\cup
\{x_0\}} (r_S-1)$. For any $S\subseteq V(W_n)$ such that $S=S_0\cup
\{x_0\}$, we have $r_S=r'_{S_0}+1$. Therefore \\$\sum_{S\subseteq
V(W_n),S=S_0\cup \{x_0\} }
\dim(\widetilde{H}_0(\Delta_S))=\sum_{S_0\subseteq V(C_n),|S_0|=i}
\dim(\widetilde{H}_0(\Delta_{S_0}))+{n\choose
i}=\\\beta_{i-1,i}(C_n) +{n\choose i}.$

The term ${n\choose i}$ is the number of subsets $S_0$ of $V(C_n)$
of cardinality $i$.

\medskip
Substituting the $\beta_{i,j}(C_n)$ from 
of \cite[Theorem~7.6.28]{ja} we have the following corollary.
\begin{cor}
Let $W_n$ be the wheel graph on $n+1$ vertices. Then the Betti numbers
of $W_n$ are as follows:\\
$(i)$ If $n=3$, then $\beta_{2,3}(W_3)=8$, $\beta_{3,4}(W_3)=3$. If $n=4$, then
$\beta_{3,4}(W_4)=9$, $\beta_{4,5}(W_4)=2$. Otherwise 
$\beta_{i,i+1}(W_n)=n{2\choose i-1}+{n\choose i}$.

\smallskip\noindent
$(ii)$  If $n=3m$, then $\beta_{2m,n}(W_n)=3m+2$, $\beta_{2m+1,n+1}(W_n)=2$.
If $n=3m+1$, then $\beta_{2m+1,n}(W_n)=3m+2$, $\beta_{2m+2,n+1}(W_n)=1$.
If $n=3m+2$, then $\beta_{2m,n}(W_n)=\beta_{2m+1,n+1}(W_n)=1$.
Otherwise, if $j>i+1$, then 
$\beta_{i,j}(W_n)=\frac{n}{n-2(j-i)}{n-2(j-i)\choose j-i}{j-i-1\choose 2i-j}$.
\end{cor}

\medskip
We can also determine the Poincar\'e series for the wheel graph algebra.
This is also a Koszul algebra, and $H_{W_n}(t)=H_{C_n}(t)+t/(1-t)$. 
Since $P_{W_n}(t)=1/H_{W_n}(-t)$ and $P_{C_n}(t)=1/H_{C_n}(-t)$, this
gives
\begin{thm}
$$P_{W_n}(t)=\frac{P_{C_n}(t)(1+t)}{1+t-tP_{C_n}(t)}$$
\end{thm}


\begin{thebibliography}{ Dillo 83}

\bibitem[B-H 98]{b-h} W.~Bruns, J.~Herzog, {\it Cohen-Macaulay rings},
revised ed., Cambridge University Press, 1998.

\bibitem[E-M-M 08]{em-mo-mo} E.~Emtander, F.~Mohammadi, S.~Moradi, 
{\it Some algebraic properties of hypergraphs}, arXiv: 0812.2366

\bibitem[F 75]{fr1} R.~Fr\"oberg, {\it Determination of a class of Poincar\'e 
series}, Math. Scand. {\bf 37}, 29--39 (1975).

\bibitem[F 78]{fr2} R.~Fr\"oberg, {\it Some complex constructions with 
applications to Poincar\'e series}, Springer Lect. Notes in Math. {\bf 740},
272--284 (1978).

\bibitem[G-L 69]{gu-le} T.~H.~Gulliksen, G.~Levin, {\it Homology of local 
rings}, Queen's paper in pure and appl. Math. {\bf 20} (1969).

\bibitem[J 04]{ja} S.~Jacques, {\it Betti Numbers of Graph Ideals}, 
Dissertation, Univ. of Sheffield (2004), arXiv:math/0410107.
\end{thebibliography}
\end{document}